\documentclass[10pt]{amsart}
\usepackage[arrow,matrix,curve,cmtip,ps]{xy}
\usepackage{amsmath,amssymb,amsfonts,amsthm,amscd}
\newtheorem{thm}{Theorem}[section]
\newtheorem{prop}[thm]{Proposition}
\newtheorem{lemma}[thm]{Lemma}
\newtheorem{coro}[thm]{Corollary}

\theoremstyle{remark}

\theoremstyle{definition}
\newtheorem{defn}[thm]{Definition}

\numberwithin{equation}{section}
\newcommand{\ip}[2]{\left\langle#1,#2  \right\rangle}

\newcommand{\norm}[1]{\left \|#1\right\|}
\newcommand{\mnorm}[1]{\left \|\left [#1\right]\right\|}
\newcommand{\abs}[1]{\left \lvert #1 \right\lvert}

\newcommand{\inj}[1]{\buildrel {\smile} \over {#1}}
\newcommand{\proj}[1]{\buildrel {\frown} \over {#1}}

\newcommand{\Min}[1]{{\rm Min}(#1)}
\newcommand{\Max}[1]{{\rm Max}(#1)}
\newcommand{\Ran}[1]{{\rm Ran}(#1)}
\newcommand{\Ker}[1]{{\rm Ker}(#1)}
\newcommand{\Ball}[1]{{\rm Ball}(#1)}

\def\proclaim #1. #2\par{\medbreak
\noindent{\bf#1.\enspace}{\sl#2}\par\medbreak}

\makeatother

\DeclareMathOperator{\htimes}{{\otimes_{\rm h}}}
\DeclareMathOperator{\A}{\mathcal{A}}
\DeclareMathOperator{\B}{\mathcal{B}}
\DeclareMathOperator{\M}{\mathcal{M}}
\DeclareMathOperator{\K}{\mathcal{K}}

\DeclareMathOperator{\Al}{{\mathcal{A}_{\ell}}}

\DeclareMathOperator{\Id}{{Id}}
\begin{document}

\title{Operator Spaces which are One-sided M-Ideals in their bidual}
\date{\today}
\author{Sonia Sharma}
\address{Department of Mathematics, University of Houston, Houston, TX
 77204-3008}
\email[Sonia Sharma]{sonia@math.uh.edu}

\begin{abstract}
We generalize an important class of Banach spaces, namely the $M$-embedded Banach spaces, to the non-commutative setting of operator spaces. 
The one-sided $M$-embedded operator spaces are the operator spaces which are one-sided $M$-ideals in their second dual.  
We show that several properties from the classical setting, like the stability 
under taking subspaces and quotients, unique extension property, Radon Nikod$\acute {\rm{y}}$m Property and many more, 
are retained in the non-commutative setting. We also discuss the dual setting of one-sided $L$-embedded operator spaces.
\end{abstract}

\maketitle

\section{Introduction}

The notion of $M$-ideals in a Banach space was introduced by Alfsen and Effros in their seminal paper
 {\cite{AE}}, in 1972. 
They defined the concept solely in terms the norm of the Banach space, deliberately avoiding any extra structure. 
Over the years the $M$-ideals have been extensively studied, resulting in a vast theory. They are an important 
tool in functional analysis. For a comprehensive treatment and for references to the extensive literature on 
the subject one may refer to the book by P. Harmand, D. Werner and W. Werner {\cite{HWW}}. Recently, the classical
 theory of $M$-ideals was generalized to the non-commutative setting of operator spaces. First, Effros and 
Ruan studied the ``complete'' $M$-ideals of operator spaces in {\cite{ER1}}. Later, Blecher, Effros, and 
Zarikian defined two varieties of $M$-ideals for the non-commutative setting, the ``left $M$-ideals'' and 
the ``right $M$-ideals'', and developed the one-sided $M$-ideal theory in a series of papers (see e.g. {\cite{BEZ, BZ1, BZ2}}
 and {\cite{BSZ}} with Smith). The intention was to create a tool for ``non-commutative functional analysis''.
 For example, one-sided $M$-ideal theory has yielded several deep general results in the theory of operator bimodules
 (see e.g. {\cite{BEZ, B3, BLM}}). In this paper, we generalize to the non-commutative setting the classical theory of an important and special 
class of $M$-ideals, called $M$-embedded spaces, namely Banach spaces which are $M$-ideals in their second dual. The study of $M$-embedded spaces marked a significant point 
in the development of $M$-ideal theory of Banach spaces. These spaces have a rich theory because of their stability 
behaviour and a natural $L$-decomposition of the third dual. These spaces have several other nice properties like 
the unique extension property, Radon Nikod$\acute{\rm{y}}$m property of the dual, 
and many more. We will study the one-sided variant of the classical theory, namely the one-sided $M$-embedded spaces. 
Our main aim is to begin to import some of the rich theory of these spaces, which comprises Chapter 3 and Chapter 4 of {\cite{HWW}}, from the classical setting to the non-commutative setting. Another motivation comes from the fact that the one-sided $M$-embedded $C^*$-algebras are nothing but 
a very nice and simple class of $C^*$-algebras, Kaplansky's ``dual $C^*$-algebras''. These are just the $C^*$-algebras of the form $\oplus^0_{i}K(H_i)$, but what is more important is that this class has many strong properties and very many interesting characterizations, see e.g. Exercise 4.7.20 from {\cite{Dix}} and references therein. We may thus try to generalize these properties and characterizations to the bigger and much more complicated
class of one-sided $M$-ideals. One-sided $M$-embedded spaces also form a subclass of the $u$-ideals of Godefroy, Kalton and Saphar {\cite{GKS}}, and so we are able to use the latter theory. We end the paper with a discussion of one-sided $L$-embedded spaces.

We will use standard operator space notations and facts. The reader is referred to the standard sources, 
see e.g. {\cite{BLM,ER2}}. We will also, frequently, use results and arguments from the one-sided $M$-ideal theory for 
which one can refer to {\cite{BEZ, BSZ, BZ2}}. 
A ({\em concrete}) operator space is a closed subspace of $B(H)$, where $H$ is some Hilbert space. Every operator 
space is completely isometrically embedded in its second dual $X^{**}$, via the canonical map $i_X: X \hookrightarrow X^{**}$.
 The canonical map $\pi_{X^*}:=i_{X^*}\circ {i_X}^*:X^{***}\longrightarrow X^{***}$ is a completely contractive projection onto 
$i_{X^*}(X^*)$ with kernel $(i_X(X))^{\perp}$. Throughout, $\inj{\otimes}$ and $\proj{\otimes}$ will denote the operator space 
injective and the operator space projective tensor products, respectively. We denote the Haagerup tensor product by $\otimes_{\rm h}$. 

Let $X\subset \A$ be an operator space, where $\A$ is a $C^*$-algebra. For each $x\in X$ 
we denote the adjoint of $x$ in $\A$ by $x^{\star}$, and the {\em adjoint} of $X$ is 
$X^{\star}=\{x^{\star} \; \mid\; x\in X\}$.
 
If $X$ is an operator space then a {\em complete left M-projection} is a linear idempotent map 
$P:X\longrightarrow X$, such that the map
\[
\nu_P^c:X\to C_2(X):x\mapsto \left [\begin{array}{c} P(x)\\ x-P(x)\end{array}\right ]
\]
 is a completely isometric injection. 
There are several equivalent characterizations of a complete left $M$-projection which can be found, 
for instance, in {\cite{BLM, BEZ, BSZ, BZ2}}. A subspace $J$ of an operator space $X$ is a {\em right M-summand} 
if it is the range of a complete left M-projection. A closed subspace $J$ of $X$ is a {\em right M-ideal} 
if $J^{\perp\perp}$ is a right $M$-summand in $X^{**}$. Similarly one may define the notion of complete right 
$M$-projections and left $M$-summands and ideals using row matrices $R_2(X)$, instead of column matrices, 
$C_2(X)$.  A right $M$-ideal which is not a right $M$-summand is called a {\em proper} right $M$-ideal. 
Mostly, we will only state the right handed version of a result since the left handed version can be proved analogously.
 Define $C_2[X]$ to be the operator space projective tensor product $C_2\proj{\otimes} X$. We say 
$P:X\longrightarrow X$ is a {\em complete right L-projection} if the map 
\[
 \nu_P^c:X\to C_2[X]: x\mapsto \left [ \begin{array}{c} P(x)\\ x-P(x)\end{array}\right ],
\]
 is a completely isometric injection. A subspace $J$ of $X$ is a left $L$-summand, if $J$ is the range of a 
complete right $L$-projection, and $J$ is a left $L$-ideal if $J^{\perp}$ is a right $M$-summand in $X^*$. Also, 
a closed subspace $J$ of $X$ is a right $M$-ideal if and only if $J^{\perp}$ is a left $L$-summand in $X^*$. Since 
a subspace $J$ of $X$ is a left $L$-ideal if and only if it is a left $L$-summand (see e.g.\ {\cite{BEZ}}), we need only 
talk about left $L$-summands.

  The {\em one-sided multipliers}, which were introduced by D. P. Blecher in {\cite{B2}}, are closely related 
to the one-sided $M$-ideals. These will be an important tool for the study of the one-sided $M$-embedded spaces.  
A projection $P$ is a complete left $M$-projection if and only if $P$ is a orthogonal projection in $\mathcal{A}_{\ell}(X)$, 
where $\mathcal{A}_{\ell}(X)$ 
is a unital $C^*$-algebra consisting of all left adjointable multipliers of $X$.
If $X$ is a dual operator space, then $\mathcal{A}_{\ell}(X)$ is a von Neumann algebra (see e.g.\ {\cite{BEZ, BLM}}).   

We now give some definitions and terminology from $u$-ideal theory which will be used in Section 3. 
If $X$ is a Banach space, then $J\subset X$ is called a {\em u-summand} if there is a 
contractive projection $P$ on $X$, mapping onto $J$, such that $\norm{I-2P}=1$. 
This norm condition is equivalent to the condition, 
\[
\norm{(I-P)(x)+P(x)}=\norm{(I-P)(x)-P(x)} \; \rm{for \; all} \; x\in X.
\]
 We call such a projection $P$, a {\em u-projection}. A subspace $J$ of $X$ is an {\em h-summand} 
if there is a contractive projection $P$ from $X$ onto $J$, such that $\norm{(I-P)-\lambda P}=1$ 
for all scalars $\lambda$ with $\abs{\lambda}=1$. This norm condition is equivalent to the condition, 
\[
\norm{(I-P)(x)-\lambda P(x)}=\norm{(I-P)(x)+P(x)} \; \rm {for \; all}\; x\in X.
\]
 Such a projection is called an {\em h-projection}. Clearly every $h$-projection is a $u$-projection
and hence every $h$-summand is a $u$-summand.

The norm condition for an $h$-summand is equivalent to saying that $P$ is hermitian in $B(X)$, that is, 
$\norm{e^{itP}}=1$ for all $t\in \mathbb{R}$. We say that $J$ is a {\em u-ideal} in $X$ if $J^{\perp}$ is a 
$u$-summand in $X^*$, and $J$ is an {\em h-ideal} if $J^{\perp}$ is an $h$-summand in $X^*$. So clearly every 
$h$-summand (resp. $u$-summand) is an $h$-ideal (resp. $u$-ideal). 
It follows from {\cite[Lemma 4.4]{BSZ}}, that one-sided $M$-summands $($ resp.\ $M$-ideals$)$ and one-sided 
$L$-summands are $h$-summands $($resp.\ $h$-ideals$)$. We refer the reader to {\cite{GKS}} for further details 
on the above topics.

\section{One-sided M-embedded spaces}
 
\begin{defn}
Let $X$ be an operator space. Then $X$ is a {\em right M-embedded} operator space 
if $X$ is a right $M$-ideal in $X^{**}$. 
We say $X$ is {\em left L-embedded} if $X$ is a left $L$-summand in $X^{**}$. 
Similarly we can define {\em right L-embedded} and {\em left M-embedded} spaces. 
If $X$ is both right and left $M$-embedded, then $X$ is called {\em completely M-embedded}. An operator space $X$ 
is {\em completely L-embedded} if $X$ is both a right and a left $L$-embedded operator space. 
\end{defn}

{\bf Remarks.} 1) \ $X$ is completely $M$-embedded if and only if $X$ is a complete $M$-ideal (in the sense of {\cite{ER1}}) in its bidual 
(see e.g. {\cite[Lemma 3.1]{BEZ}} and {\cite[Chapter 7]{BZ2}}).
 
2) \ Reflexive spaces are automatically completely $M$-embedded. Let $X$ be a right M-summand in $X^{**}$. 
Since  $X^{**}$ is a dual operator space, $X$ is $w^*$-closed, see {\cite[p.8]{BZ2}}.
 So $X = X^{**}$. Hence a non-reflexive operator space cannot be a non-trivial one-sided $M$-summand in its second dual. 
So non-reflexive  right $M$-embedded spaces are proper right $M$-ideals.  
Henceforth we will assume all our operator spaces to be non-reflexive.

\medskip

We state an observation of David Blecher which provides an alternative definition of completely $L$-embedded operator 
spaces. To explain the notation here, $M_n(X^*)_*$ is the `obvious' predual
of $M_n(X^*)$, namely the operator space
projective tensor product of the predual of $M_n$ and $X$.

\begin{lemma}\label{Lproj}
Let $X$ be an operator space. Then there exists a complete $L$-projection from $X^{**}$ onto $X$ 
if and only if for each $n$, there exists a $L$-projection from $M_n(X^*)^*$ onto $M_n(X^*)_*$.
\end{lemma}
\begin{proof}
We are going to use the well known principle that if $J$ is a subspace
of $X$, then $J$ is an $L$-summand (resp.\ left $L$-summand, complete
$L$-summand) of $X$ iff $J^\perp$ is an $M$-summand (resp.\ right $M$-summand, complete
$M$-summand) of $X^*$. See for example the proof of \cite[Proposition 3.9]{BEZ}.

By the above principle, $X$ is a complete $L$-summand of $X^{**}$
iff $X^\perp$ is a complete $M$-summand of $X^{***}$.
By {\cite[Proposition 4.4]{ER1}}, this happens iff $M_n(X^\perp)$ is an
$M$-summand of $M_n(X^{***})$ for each $n$. Now
$M_n(X^{***})$ is the dual of the operator space
projective tensor product of the predual of $M_n$ and $X^{**}$. 
Moreover, $M_n(X^\perp)$ is easily seen to be the `perp' of
the operator space
projective tensor product of the predual of $M_n$ and $i_X(X)$.
That is, $M_n(X^\perp) = (M_n(X^*)_*)^\perp$. (We are using facts from {\cite[Proposition 7.1.6]{ER2}}). 
By the above principle,
we deduce that $X$ is a complete $L$-summand of $X^{**}$ iff
$M_n(X^*)_*$ is a $L$-summand of $M_n(X^*)^*$ for each $n$.
\end{proof}

\begin{prop}\label{ns}
Let $X$ be an operator space, then the following are equivalent: 
\begin{enumerate} 
\item[(i)] $X$ is a right $M$-ideal in $X^{**}$.
\item[(ii)] The natural projection $\pi_{X^*}$ is a complete right $L$-projection.
\end{enumerate}
\end{prop}
\begin{proof}
(i)$\Rightarrow$(ii) \ Let $X\cong i_X(X)$ be a right $M$-ideal in $X^{**}$, 
then $i_X(X)^{\perp}$ is a complete left $L$-ideal in $X^{***}$. Let $P$ be a complete right $L$-projection 
onto $i_X(X)^{\perp}$, then  $i_X(X)^{\perp}$ is the kernel of the complementary right $L$-projection, namely
 $I-P$. Now  Ker $\pi_{X^*}=(i_X(X))^{\perp}=$ Ker$(I-P)$. So by {\cite [Theorem 3.10(b)]{BEZ}}, $\pi_{X^*}=I-P$. 
Hence $\pi_{X^*}$ is a complete right $L$-projection.

(ii)$\Rightarrow$(i) \ If  $\pi_{X^*}$ is a complete right $L$-projection, then so is $I-\pi_{X^*}$. Now
\[
\Ran {I-\pi_{X^*}}=\Ker {\pi_{X^*}}=(i_X(X))^{\perp}.
\]
 So $(i_X(X))^{\perp}$ is a left $L$-summand in $X^{***}$, and hence $i_X(X)$ is a right $M$-ideal in $X^{**}$.   
\end{proof}

\begin{coro}\label{coro_duality}
 If $X$ is a right $M$-embedded operator space, then $X^*$ is a left $L$-embedded operator space.
\end{coro}
\begin{proof}
 Since $i_{X^*}(X^*)$ is the range of $\pi_{X^*}$, and by Proposition {\ref{ns}}, $\pi_{X^*}$ is a complete right $L$-projection on $X^{***}$, the result follows. 
\end{proof}

However, it is not true that if $X$ is a left $L$-summmand in its bidual then $X^*$ is a right $M$-summand in its bidual. For example take $X=S^1(H)$, the trace class operators on $H$. Then since $\K(H)$ is a complete $M$-ideal in $B(H)$, by Corollary {\ref{coro_duality}}, $S^1(H)$ is complete $L$-summand in $B(H)^*$. But $B(H)$ is not a right (or left) $M$-summand in $B(H)^{**}$ since $B(H)$ is non-reflexive.

\begin{prop}
If $X$ is a $M$-embedded Banach space, then $\Min X$ is a completely $M$-embedded operator space. If $X$ is $L$-embedded, then $\Max X$ is completely $L$-embedded.
\end{prop}
\begin{proof}
Let $X$ be a $M$-ideal in $X^{**}$, then $\Min X$ is a two-sided $M$-ideal in $\Min {X^{**}}$. Indeed if $Z$ is a Banach space, then the right $M$-ideals, as well as the left $M$-ideals, of $\Min Z$, coincide with the $M$-ideals of $Z$ (see e.g.\ {\cite{BEZ}}). But $\Min {X^{**}}=\Min X^{**}$ completely isometrically. So $\Min X$ is a right $M$-ideal in $\Min X^{**}$, and hence $\Min X$ is $M$-embedded. The second assertion follows similarly, using the fact that $L$-ideals of any Banach space $Z$ coincide with the left, as well as the right, $L$-ideals of $\Max Z$.
\end{proof}

\begin{thm}{\label{stability}}
Let $X$ be a right $M$-embedded space and $Y$ be a subspace of $X$, then both $Y$ and $X/Y$ are right $M$-embedded.
\end{thm}
\begin {proof}
We first show that $Y$ is right $M$-embedded. By Proposition \ref{ns}, we need to show that $\pi_{Y^*}$ is a complete right $L$-projection. Let $i:Y\longrightarrow X$ be the inclusion map, then $i^{***}$ is a complete quotient map. So for every $[v_{i j}]\in M_n(Y^{***})$ we can find $[w_{ij}]\in M_n(X^{***})$ such that, $i^{***}_n([w_{ij}])=[v_{ij}]$ and $\norm{[w_{ij}]}\leq \norm{[v_{ij}]}$. Also note that $\pi_{Y^*}\circ i^{***}=i^{***}\circ \pi_{X^*}$. For $[v_{ij}]$ and $[w_{ij}]$  as above, we have
\vspace{.15cm}

$\begin{array}{rll}
\hspace{1.95cm}&& \mnorm{\pi_{Y^*}(v_{ij})\;\;\;\;\;\;v_{ij}-\pi_{Y^*}(v_{ij})}_{M_n({R_2[Y^{***}]})}\\
&=&\mnorm{\pi_{Y^*}i^{***}(w_{ij})\;\;\;\;\;\;i^{***}(w_{ij})-\pi_{Y^*}i^{***}(w_{ij})}_{M_n(R_2[Y^{***}])}\\
&=&\mnorm{i^{***}\pi_{X^*}(w_{ij})\;\;\;\;\;\;i^{***}(w_{ij})-i^{***}\pi_{X^*}(w_{ij})}_{M_n({R_2[Y^{***}]})}\\
&\leq & \norm{i^{***}}_{cb}\mnorm{\pi_{X^*}(w_{ij})\;\;\;\;\;\;w_{ij}-\pi_{X^*}(w_{ij})}_{M_n({R_2[X^{***}]})}\\
&=&\mnorm{w_{ij}}\\
&\leq & \mnorm{v_{ij}}.
\end{array}\\$

 This shows that the map $\mu_{\pi_{Y^*}}^{r}:Y^{***}\longrightarrow R_2[Y^{***}]$ given by
\[
\mu_{\pi_{Y^*}}^{r}(y)=[\pi_{Y^*}(y)\;\;\;\;\;\; y-\pi_{Y^*}(y)],
\]
is a complete contraction.
 Now since $\proj{\otimes}$ is projective, and $i^{***}$ is a complete quotient map, then so is $i^{***}\otimes Id:R_2[X^{***}]\longrightarrow R_2[Y^{***}]$. For each $[\;y_{ij} \;\;\;\;\;\acute{y_{ij}}\; ]\in  R_2[Y^{***}]$ we can find $[x_{ij}\;\;\;\;\;\acute{ x_{ij}}]\in R_2[X^{***}]$, such that  $(i^{***}\otimes Id)([\;x_{ij}\;\;\;\;\; \acute{x_{ij}}\;])=[y_{ij} \;\;\;\;\;\acute{y_{ij}}]$ and $\norm{[x_{ij} \;\;\;\;\;\acute{x_{ij}} ]}\leq \norm{[\;y_{ij}\;\;\;\;\; \acute{y_{ij}}\;]}$. Consider\\

$\begin{array}{lll}
\hspace{1.5cm}&&\mnorm{\pi_{Y^*}(y_{ij})+\acute{y_{ij}}-\pi_{Y^*}(\acute{y_{ij}})}_{M_n({[Y^{***}]})}\\
&=&\mnorm{\pi_{Y^*}(i^{***}(x_{ij}))+i^{***}(\acute{x_{ij}})-\pi_{Y^*}(i^{***}(\acute{x_{ij}}))}_{M_n({[Y^{***}]})}\\
&=&\mnorm{i^{***}\pi_{X^*}(x_{ij})+i^{***}(\acute{x_{ij}})-i^{***}\pi_{X^*}(\acute{x_{ij}})}_{M_n({[Y^{***}]})}\\
&\leq&\mnorm{\pi_{X^*}(x_{ij})+(\acute{x_{ij}})-\pi_{X^*}(\acute{x_{ij}})}_{M_n({[X^{***}]})}\\
&\leq& \mnorm{x_{ij} \;\;\;\;\;\;\acute{x_{ij}}}_{M_n({R_2[X^{***}]})}\\
&\leq& \mnorm{y_{ij} \;\;\;\;\;\;\acute{y_{ij}}}_{M_n({R_2[Y^{***}]})}.
\end{array}\\$
\vspace{.25cm}

 This shows that the map $\nu_{\pi_{Y^*}}^r:R_2[Y^{***}]\longrightarrow Y^{***}$ given by $$\nu_{\pi_{Y^*}}^r([\;y \;\;\;\; \acute{y}\;])=\pi_{Y^*}(y)+\acute{y}-\pi_{Y^*}(\acute{y}),$$ is a complete contraction. Hence by {\cite[Proposition 3.4]{BEZ}}, $\pi_{Y^*}$ is a complete right $L$-projection.

Consider the canonical complete quotient map $q:X\longrightarrow X/Y$, then $q^{***}: (X/Y)^{***}\longrightarrow X^{***}$ is a complete isometry. We also have that $\pi_{(X/Y)^*}\circ q^{***}=q^{***}\circ \pi_{(X/Y)^*}$. Since $R_2[(X/Y)^{***}]=R_2\otimes_{\rm h} (X/Y)^{***}$ and $R_2[X^{***}]=R_2\otimes_{\rm h} X^{***}$, and $\otimes_{\rm h}$ is injective, the map $Id\otimes q^{***}:R_2[(X/Y)^{***}]\longrightarrow R_2[X^{***}]$ will be a complete isometry. We need to show that $\pi_{(X/Y)^*}$ is a complete right $L$-projection on $(X/Y)^{***}$. For the sake of convenience we will write $\pi$ for $\pi_{(X/Y)^*}$. Let $[v_{i j}]\in M_n((X/Y)^{***})$, then by using the above facts we get\\ 

$\begin{array}{lll}
\hspace{1.5cm}&&\mnorm{\pi(v_{i j})\;\;\;\;\;\;v_{i j}-\pi(v_{i j})}_{M_n(R_2[(X/Y)^{***}])}\\
&=&\mnorm{(q^{***}\circ \pi)(v_{i j})\;\;\;\;\;\;q^{***}(v_{i j})-(q^{***}\circ \pi)(v_{i j})}_{M_n(R_2[X^{***}])}\\
&=&\mnorm{(\pi\circ q^{***})(v_{i j})\;\;\;\;\;\;q^{***}(v_{i j})-(\pi\circ q^{***})(v_{i j})}_{M_n(R_2[X^{***}])}\\
&=& \mnorm{q^{***}(v_{i j})}_{M_n(X^{***})}\\
&=& \mnorm{v_{i j}}_{M_n((X/Y)^{***})}.
\end{array}\\$

This shows that $\pi_{(X/Y)^*}$ is a complete right $L$-projection. Since Ran($\pi_{(X/Y)^*})=(X/Y)^*$, $X/Y$ is right $M$-embedded.
\end{proof}

{\bf Remark.} The property of one-sided ``$M$-embeddedness'' of subspaces and quotients does not pass to extensions, i.e., if $Y$ is a subspace of $X$ such that $Y$ and $X/Y$ are right $M$-embedded spaces, then $X$ need not be right $M$-embedded. Consider $X=c_0\oplus_1 c_0$ and $Y=c_0\times \{ 0\}$, both with minimal operator space structure. Since $Y$ and $X/Y$ are $M$-embedded, $\Min {Y}$ and $\Min {X/Y}$ are completely $M$-embedded.  
Let $P$ be the contractive projection from $X$ onto $Y$, then $I-P$ is completely contractive, and hence a complete quotient map, from $\Min {Y}$ onto 
$\Min {\Ran {I-P}}$. Thus $\Min Y /\Ker P \cong \Min {\Ran {I-P}}$, 
completely isometrically. But $\Ran {I-P}= Y/X$ isometrically, so $\Min{X}/\Min{Y}\cong \Min {X/Y}$, completely isometrically. 
Now if $\Min{X}^{**}$ has a non-trivial right $M$-ideal, then since $\Min{X}^{**}=\Min{X^{**}}$, $X^{**}$ has a nontrivial $M$-ideal. But this is not possible, since $X^{**}$ has a non-trivial $L$-summand, and by {\cite [Theorem I.1.8]{HWW}}, a Banach space cannot contain nontrivial $M$-ideals and nontrivial $L$-summands simultaneously, unless it is two dimensional.

\begin{prop}{\label{matrix_const}}
Let $X$ be a left $($right$)$ $M$-embedded space, then
\begin{enumerate} 
\item[(i)] $M_{m,n}(X)$ is left $($right$)$ $M$-embedded for all $m$ and $n$. In particular, $C_n(X)$ $($resp.\ $R_n(X)$$)$ is left $($right$)$ $M$-embedded in $C_{n}(X^{**})$ $($resp.\ $R_n(X^{**}))$. 
\item [(ii)]  $C_{\infty}(X)$ $($resp.\ $R_{\infty}(X)$$)$ is a left $($right$)$ $M$-ideal in $C_{\infty}(X^{**})$ $($resp.\ $R_{\infty}(X^{**}))$.
\end{enumerate}
\end{prop}
\begin{proof}
(i) \ If $J\subset X$ is a right $M$-ideal then $M_{m,n}(J)$ is a right $M$-ideal in $M_{m,n}(X)$ (see e.g.\ {\cite{BEZ}}). Now the result follows from the fact that $M_{m,n}(X^{**})$
$=M_{m,n}(X)^{**}$ completely isometrically.
    
(ii) \ If $X$ is a left $M$-ideal in $X^{**}$, then by the left-handed version of Theorem 5.38 from {\cite{BZ2}}, $C_{\infty}\otimes_{\rm h} X$ is a left $M$-ideal in $C_{\infty}\otimes_{\rm h} X^{**}$. But  $C_{\infty}\otimes_{\rm h} X=C_{\infty}\inj{\otimes} X=C_{\infty}(X)$ and  $C_{\infty}\otimes_{\rm h} X^{**}=C_{\infty}\inj{\otimes} X^{**}=C_{\infty}(X^{**})$. For the second assertion, use {\cite [Theorem 5.38]{BZ2}} and that $Y\otimes_{\rm h} R_{\infty}=R_{\infty}(Y)$, for any operator space $Y$.
\end{proof}

It would be interesting to know when is $C_{\infty}(X)$ a right $M$-embedded space, that is, whether one can replace $C_{\infty}(X^{**})$ by $C_{\infty}(X)^{**}=C_{\infty}^w(X^{**})$ in Proposition {\ref{matrix_const}} (ii)? We will see in the remark after Proposition {\ref{TROresult}} that this is true in case of TROs. Also note that, if $X$ is a WTRO then a routine argument shows that $C_{\infty}(X)$ is a right $M$-ideal in $C_{\infty}^{w}(X)$.  

\begin{prop}
 Every right $M$-embedded $C^*$-algebra is left $M$-embedded. 
\end{prop}
\begin{proof}
Suppose $\A$ is a right $M$-ideal in $\A^{**}$, then $\A$ is a closed right ideal in $\A^{**}$ and $\A^{\star}$ is a closed left ideal in $\A^{**}$. But $\A$ is self-adjoint, i.e., $\A=\A^{\star}$, hence $\A$ is a two-sided $M$-ideal in $\A^{**}$ (see e.g.\ {\cite[Section 4.4]{BZ2}}).
\end{proof}

{\bf Remark.} A complete $M$-ideal in an operator space is an $M$-ideal in the underlying Banach space. So by the above proposition, a one-sided $M$-embedded $C^*$-algebra is a $M$-embedded $C^*$-algebra in the classical sense. Hence by {\cite[Proposition III.2.9]{HWW}}, it has to be $*$-isomorphic to  $\oplus_i^{0}\mathcal{K}(H_i)$ (a $c_{0 }$-sum), for Hilbert spaces $H_i$. These are Kaplansky's dual $C^*$-algebras, consequently, one-sided $M$-embedded $C^*$-algebras satisfy a long list of equivalent conditions which can be found for instance in the works of Dixmier and Kaplansky (see e.g.\ Exercise 4.7.20 from {\cite{Dix}}). To mention a few:

\begin{enumerate}
\item[(i)] Every closed right ideal $J$ in $\A$ is of the form $e\A$ for a projection $e$ in the multiplier algebra of $\A$.
\item[(ii)] There is a faithful $*$-representation $\pi:\A\longrightarrow \mathcal{K}(H)$ as compact operators on some Hilbert space $H$.
\item[(iii)] The sum of all minimal right ideals in $\A$ is dense in $\A$.
\end{enumerate} 

We imagine that several of these have variants that are valid for general one-sided $M$-embedded spaces (see e.g.\ Theorem {\ref{thm_centralizer}}).

\medskip

Note that in contrast to the fact that every one-sided $M$-embedded $C^*$-algebra is completely $M$-embedded, we give examples in {\cite{ABS}} of non-selfadjoint left $M$-embedded algebras which are not right $M$-embedded.

\begin{thm}{\label{example}}
Let $Y$ be a non-reflexive operator space which is right $M$-embedded and if $X$ is any finite dimensional operator space, then $Y\otimes_{\rm h} X$ is right $M$-embedded. Further, if $Z(Y^{(4)}\otimes_{\rm h}X) \cong \mathbb{C}I$ then $Y\otimes_{\rm h} X$ is not left $M$-embedded, where $Z(X)$ denotes the centralizer algebra of $X$.
\end{thm} 
\begin{proof}
Since $Y$ is a right $M$-ideal in $Y^{**}$, by {\cite [Proposition 5.38]{BZ2}}, $Y\otimes_{\rm h} X$ is a right $M$-ideal in $Y^{**}\otimes_{\rm h} X$. Since $X$ is finite dimensional, $(Y\otimes_{\rm h} X)^{**}=Y^{**}\otimes_{\rm h}X$ (see e.g.\ {\cite[1.5.9]{BLM}}). Hence $Y\otimes_{\rm h} X$ is a right $M$-ideal in its bidual. 
Suppose that $Y\otimes_{\rm h} X$ is also left $M$-embedded and $P$ be a projection in $ Z(Y^{(4)}\otimes_{\rm h} X)$ such that $(Y\otimes_{\rm h} X)^{\perp\perp}=P(Y^{(4)}\otimes_{\rm h} X)$. Since $Z(Y^{(4)}\otimes_{\rm h}X) \cong \mathbb{C}I$, so $(Y\otimes_{\rm h} X)^{\perp\perp}=Y^{(4)}\otimes_{\rm h} X$. 
Now note that for any operator space $X$, if $E=i_{X}(X)\subset X^{**}$, then by basic functional analysis, $i_{X^{**}}(X^{**})\cap E^{\perp\perp}=i_{X^{**}}(E)$. So if $i_{X}(X)^{\perp\perp}=X^{(4)}$, then $i_{X^{**}}(E)=i_{X^{**}}(X^{**})$, hence $X^{**}=E=i_{X}(X)$. 
This implies that $Y^{**}\otimes_{\rm h} X\cong Y\otimes X$, which is not possible since $Y$ is non-reflexive.         
\end{proof}

 As a result, we can generate many more concrete examples of right $M$-embedded spaces which are not left $M$-embedded. If $\A$ is any algebra of compact operators, e.g.\ a nest algebra of compact operators, then we know that it is two-sided $M$-embedded. Hence, $Z=\A\otimes_{\rm h} X$ is right $M$-embedded for all finite dimensional operator spaces $X$, as are all subspaces of $Z$. Almost all of these will, surely, not be left $M$-embedded. We show in {\cite{ABS}} that if $\A$ and $\B$ are approximately unital operator algebras, then $ \Al (\A \htimes \B)\cong \Delta(M(\A))$, if $\Delta(\A)$ is not one dimensional. As a consequence, using a similar argument as in Theorem {\ref{example}}, we can show that if $\A$ and $\B$ are approximately unital operator algebras such that $\A$ is completely $M$-embedded and $\B$ is finite dimensional with $\B\ne \mathbb{C}1$, then $\A\htimes \B$ is a right $M$-embedded operator space which is not left $M$-embedded. 
 
\medskip

For an operator space $X$, the {\em density character} of $X$ is the least cardinal $m$ such that there exists a dense subset $Y$ of $X$ with cardinality $m$. We denote the density character by {\em dens$(X)$}. So if $X$ is separable, then {\em dens$(X)=\aleph_{\circ}$}.  A Banach space $X$ is an {\em Asplund space} if every separable subspace has a separable dual. Also $X$ is an Asplund space if and only if $X^*$ has the RNP. For more details see {\cite[p.91, p.132]{BOU}}, {\cite[p.82, p.195, p.213]{DU}} and {\cite[p.34, p.75]{Phe}}.
Using an identical argument to the classical case (see {\cite[Theorem III.3.1]{HWW}}), we can show the following. Note that the following proof uses a couple of ideas from Section 3 such as the notion of a norming subspace and Corollary {\ref{norming}}.  
\begin{thm}{\label{RNP}}
If $X$ is right $M$-embedded and $Y$ is a subspace of $X$, 
then dens$(Y)$ $=$ dens$(Y^*)$. In particular, separable subspaces of $X$ have a separable dual. So right $M$-embedded spaces are Asplund spaces, and $X^*$ has the Radon-Nikod$\acute{ y}$m Property.  
\end{thm}

\begin{proof}
By Proposition {\ref {stability}}, we can assume WLOG that $Y=X$. Suppose that $K$ is a dense subset of $X$. For each $x\in K$ choose $x^*\in X^*$ such that $\norm{x^*}=1$ and $x^*(x)=\norm{x}$. Then $N=\overline{\rm {span}}\{x^*\; :\;x\in K \}$ is norming for $X^*$, but by Corollary {\ref{norming}}, $X^*$ has no nontrivial norming subsets. So $N=X^*$ and hence dens$(X^*)=$dens$(X)$.    
\end{proof}

A TRO is a closed subspace $X$ of a $C^*$-algebra such that $XX^{\star}X\subset X$. A WTRO is a $w^*$-closed subspace of a von Neumann algebra with $XX^{\star}X\subset X$. A TRO is essentially the same as a 
Hilbert $C^*$-module (see e.g.\ {\cite[8.1.19]{BLM}}). If $X$ is a TRO, then $X$ is a Hilbert $C^*$-bimodule 
over $XX^{\star}$-$X^{\star}X$ (see e.g.\ {\cite[8.1.2]{BLM}}).

\begin{lemma}{\label{lemmaTROci1}}
If $Z$ is a TRO which is isometrically isomorphic to $\mathcal{K} (H,K)$, then $Z$ is completely isometrically isomorphic to either $\mathcal{K}(H,K)$ or $\mathcal{K}(K,H)$. 
\end{lemma}
\begin{proof}
Let $\theta: Z \longrightarrow  \mathcal{K} (H,K)$ be an isometric isomorphism. Then $\theta^{**}: Z^{**}\longrightarrow B(H,K)$ is an isometric isomorphism. We use {\cite[Theorem 2.1]{Sol}} to prove this result. 
Take $M=B(K\oplus H)$ and $N=\mathcal{L}(Z)^{**}$, where $\mathcal{L}(Z)$ denotes the linking $C^*$-algebra of $Z$. Then by Lemma 3.1 from {\cite{Sol}}, there exists a projection $q\in \mathcal{L}(Z)^{**}$ such that both $q$ and $I-q$ have central support equal to $I$, and $q N(I-q) \cong Z^{**}$. Let $p=P_K\in B(K)$, then $pM(I-p)\cong B(H,K)$.
Thus by {\cite[Theorem 2.1]{Sol}}, there exist central projections $e_1$, $e_2$ in $M$ and $f_1$, $f_2$ in $N$ with $e_1 +e_2= I_{K\oplus H}$ and $f_1+ f_2=I_{\mathcal{L}(Z)^{**}}$. 
Since there are no central projections in $B(K\oplus H)$, either $e_1=I_{K\oplus H}$ or $e_1=0$. By {\cite[Theorem 2.1]{Sol}}, either there exists a $*$-isomorphism 
$\psi: B(K\oplus H)\longrightarrow f_1 N f_1$ such that $(\theta^{**})^{-1}=\psi|_{B(H,K)}$, or there exists a $*$-`anti'-isomorphism $\phi: B(K\oplus H)\longrightarrow f_2 N f_2$ such that $(\theta^{**})^{-1}=\phi|_{B(H,K)}$. In the first case, $\psi$ is a complete isometry and hence so is $(\theta^{**})^{-1}$. Thus $Z$ is completely isometrically isomorphic to $\mathcal{K}(H,K)$. We claim that the second case implies that $Z$ is completely isometrically isomorphic to $\mathcal{K}(K,H)$. Let $\{e_i\}$ and $\{f_j\}$ be orthonormal bases for $K$ and $H$ resp., then $S=\{e_i\}\ \cup \{f_j\}$ is an orthonormal basis for $K\oplus H$. For each $T\in B(K\oplus H)$, define $\tilde{T}\in B(K\oplus H)^{\rm{op}}$ to be the transpose of $T$ given by $\tilde{T}\eta=\sum_i \ip{Te_i}{\eta}e_i+\sum_j \ip{Tf_j}{\eta}f_j$, for every $\eta\in S$. Then $t: B(K\oplus H)\longrightarrow B(K\oplus H)^{\rm{op}}$ defined as $t(T)=\tilde{T}$, is a $*$-`anti'-isomorphism and $t(B(K,H))=B(H,K)$. So $\tilde{\phi}=\phi\circ t$ is a $*$-isomorphism, and hence a complete isometry, such that $\tilde{\phi}(B(K,H))=\phi(B(H,K))=(\theta^{**})^{-1}(B(H,K))$. Thus restriction of $\tilde{\phi}$ to $B(K,H)$ is a complete isometry onto $Z^{**}$. 
\end{proof}

\begin{prop}{\label{TROresult}}
A one-sided $M$-embedded TRO is completely isometrically isomorphic to the $c_0$-sum of compact operators on some Hilbert spaces.
\end{prop}
\begin{proof}
Let $X$ be a right $M$-embedded TRO, then by Theorem {\ref{RNP}}, $X^*$ has the RNP. Also since $X$ is a TRO, it is a $JB^*$-triple. From {\cite{BG}} we know that if $X$ is a $JB^*$-triple and $X^*$ has the Radon-Nykod$\acute{\rm y}$m property, then $X^{**}$ is isometrically an $l^{\infty}$-sum of type-I triple factors, i.e., $X^{**}\cong \oplus^{\infty}_i B(H_i,K_i)$ isometrically, for some Hilbert spaces $H_i$ and $K_i$.  
By Proposition {\ref{bitranspose}} (ii), there exists a surjective isometry $\rho: X\longrightarrow \oplus^{0}_i\mathcal{K} (H_i,K_i)$.    
Let $\K_i=\mathcal{K}(H_i,K_i)$, $\rho_i=\rho^{-1}|_{\K_i}$ and $Z_i=\rho_i(\K_i)$, then $X\cong \oplus_i^0 Z_i $, isometrically. 
So each $Z_i$ is a $M$-summand in $X$. Every $M$-summand in a TRO is a complete $M$-summand, and hence each $Z_i$ is a sub-TRO of $X$ (see e.g.\ {\cite[8.5.20]{BLM}}). Also since the $Z_i$ are orthogonal, there is a ternary isomorphism between $ \oplus_i^0 Z_i$ and $X$, given by $(x_i)\mapsto \sum_i x_i$. Hence $X\cong \oplus_i^0 Z_i $ completely isometrically (see e.g.\ {\cite[Lemma 8.3.2]{BLM}}).
 Thus by Lemma \ref{lemmaTROci1}, for each $i$, either $\rho^{-1}_i$ is a complete isometry or there exists a complete isometry $\tilde{\rho}_i:Z_i\longrightarrow \mathcal{K}(K_i,H_i)$. 
Define $\theta=\oplus \pi_i: \oplus^{\infty}_iZ_i^{**}\longrightarrow \oplus_i^{\infty} B_i$, where $\pi_i$ is either $(\rho_i^{-1})^{**}$ or $(\tilde{\rho}_i)^{**}$ and $B_i$ is either $B(H_i,K_i)$ or $B(K_i,H_i)$.  
Since each $Z_i$ is a WTRO, and each $\pi_i$ a complete isometry,  
$\theta$ is a complete isometry. So $X^{**}$ is completely isometrically isomorphic to $\oplus_i^{\infty} B_i$. Hence by Proposition \ref{bitranspose} (ii), $X$ is completely isometrically isomorphic to $\oplus_i^{0}\mathcal{K}_i$ where $\mathcal{K}_i$ is either $\mathcal{K}(H_i,K_i)$ or $\mathcal{K}(K_i,H_i)$.
\end{proof}

{\bf Remark.} As we stated earlier, if $X$ is a right $M$-embedded TRO then $C_{\infty}(X)$ is also right $M$-embedded. Indeed, by Proposition {\ref{TROresult}}, $X$ is completely isometrically isomorphic to $\oplus_i^0\mathcal{K}(H_i,K_i)$, which implies that $C_{\infty}(X)$ is completely isometrically isomorphic to $\oplus_i^0\mathcal{K}(H_i,K_i^{\infty})$, and the latter is a complete $M$-ideal in its bidual. More generally, if $X$ is a right $M$-embedded TRO, then every $\mathbb{K}_{I,J}(X)$ is completely $M$-embedded.

\section{More Properties of one-sided $M$-embedded spaces}

In this section we show that one-sided $M$-embedded spaces retain a number of properties of the classical $M$-embedded spaces. 
We start by stating two theorems from $u$-ideal theory which will allow us to draw some useful conclusions about one-sided $M$-embedded spaces. For proof of these theorems see {\cite[Theorem 6.6]{GKS}} and {\cite [Theorem 5.7]{GKS}}, respectively. A $u$-ideal $J$ of $X$ is a strict $u$-ideal if $\Ker{P}$ is a norming subspace in $X^*$, where $P$ is a $u$-projection onto $J^{\perp}$. By a {\em norming subspace} of $X^*$ we mean a subspace $N$ of $X^*$ such that for each $x\in X$, $\norm{x}=\sup \{ \abs{\phi(x)}: \phi\in N,\; \norm{\phi}\leq 1 \}$.

\begin{thm}{\label{godf1}}
 Let $X$ be a Banach space which is an h-ideal in its bidual. Then the following are equivalent:\\
 $($a$)$ $X$ is a strict h-ideal. $($b$)$ $X^*$ is an h-ideal. $($c$)$ $X$ contains no copy of $\ell^1$.
\end{thm}
\begin{thm}{\label{godf2}}  
 Let $X$ be a Banach space which contains no copy of $\ell^1$ and is a strict u-ideal, then 
\begin{enumerate}
\item [(a)] if $T:X^{**}\longrightarrow X^{**}$ is a surjective isometry, then $T=S^{**}$, for some surjective isometry  $S:X\longrightarrow X$. 
\item [(b)] $X$ is the unique isometric predual of $X^*$ which is a strict u-ideal.
\end{enumerate}
\end{thm}

\begin{prop}\label{bitranspose}
Let X be a non-reflexive right M-ideal in its bidual, then 
\begin{enumerate}
\item[(i)] $X$ do not contain a copy of $\ell^1$ and $X$ is a strict $u$-ideal.
\item[(ii)] If $T:X^{**}\longrightarrow X^{**}$ is a $($completely$)$ isometric surjection, then $T$ is a bitranspose of some $($completely$)$ isometric surjective map on $X$. 
\item[(iii)] $X$ is the unique isometric predual of $X^*$.
\end{enumerate}
\end{prop}
\begin{proof}
(i) \ Suppose that $X$ is a right $M$-ideal in $X^{**}$, then being the range of a complete right $L$-projection, $X^{*}$ is a right $L$-summand. So $X$ and $X^*$ are both $h$-ideals (see Section 1). Hence, by Theorem {\ref{godf1}}, $X$ is a strict $u$-ideal and $X$ does not contain a copy of $\ell^1$. 

(ii) \ By Theorem {\ref{godf2}} and (i), $T=S^{**}$ for some isometric surjection $S$ on $X$. Further, if $T$ is a complete isometry then it is not difficult to see that $S$ is also a complete isometry.
 
(iii) \ This follows from Theorem {\ref{godf2}} and (i).    
\end{proof}

\begin{thm}\label{thm_centralizer}
Suppose that $X$ is a left $M$-embedded operator space. Then
\begin{itemize} 
\item[(i)] Every right $M$-ideal of $X$ is a right $M$-summand.
\item[(ii)] Every complete left $M$-projection $P$ in $X^{**}$ is the bitranspose of a complete left $M$-projection $Q$ on $X$.
\item[(iii)] Suppose $X$ is also a right $M$-embedded operator space and it has no nontrivial right $M$-summands. Then for every nontrivial right $M$-ideal $J$ of $X^{**}$, either $J$ contains $X$ or $J\cap X=\{0\}$. 
\item[(iv)] $\mathcal{A}_{\ell}(X)\cong \mathcal{A}_{\ell}(X^{**})$. 
\item[(v)] If $X$ is a completely $M$-embedded space, then ${Z}(X)={Z}(X^{**})$, where ${Z}(X)$ is the centralizer algebra of $X$, in the sense of {\cite [Chapter 7]{BZ2}}.
\end{itemize}
\end{thm}
\begin{proof} 
(i) \ Let $J$ be a right $M$-ideal of $X$ and suppose that $P$ is a projection in Ball$(\mathcal{M}_{\ell}(X^{**}))$ such that $P(X^{**})=J^{\perp\perp}$. Since $X$ is a left $M$-ideal in $X^{**}$ and $P\in \mathcal{M}_{\ell}(X^{**})$, by \cite[Proposition 4.8]{BZ2} we have that $P(X)\subset X$. Then $Q:=P|_X \in \mathcal{M}_{\ell}(X)$ and $\norm{Q}=\norm{P|_X}_{\mathcal{M}_{\ell}(X)}\leq 1$. Also $J^{\perp\perp}\cap X=J$ is the range of $Q$. Hence by {\cite[Theorem 5.1]{BEZ}}, $J$ is a complete right $M$-summand.
 
(ii) \ If $P$ is an complete left $M$-projection in $X^{**}$, then $T:=2P-Id_{X^{**}}$ is a complete isometric surjection of $X^{**}$. Hence by Proposition \ref{bitranspose}, $T=S^{**}$, for some complete surjective isometry on $X$. So, $2P=T+Id_{X^{**}}=(S+Id_X)^{**}$, and since by {\cite[Section 5.3]{BZ2}} $\mathcal{A}_{\ell}(X)\subset \mathcal{A}_{\ell}(X^{**})$, $S+Id_{X}$ must be a complete left $M$-projection in $X$.

(iii) \ Let $P$ be a complete two-sided $M$-projection from $X^{(4)}$ onto $X^{\perp\perp}$ and $Q$ be a complete left $M$-projection from $X^{(4)}$ onto $J^{\perp\perp}$. Then by \cite[Theorem 5.1]{BEZ}, $P\in$ Ball($\mathcal{M}_{\ell}(X^{(4)})$) and $Q\in$ Ball($\mathcal{M}_r(X^{(4)})$), which implies that $PQ=QP$. Hence by {\cite[Theorem 5.30 (ii)]{BZ2}}, $J\cap X$ is a right $M$-ideal in $X^{**}$. But $J\cap X \subset X$, so by {\cite[Theorem 5.3]{BZ2}}, $J\cap X$ is a right $M$-ideal in $X$. Hence by (i), $J\cap X$ is a right $M$-summand. By the hypothesis, either $J\cap X=\{ 0 \}$ or $J\cap X=X$, i.e., $J\cap X=\{0\}$ or $X\subset J$.   

(iv) \ We know that $\mathcal{A}_{\ell}(X)\subset \mathcal{A}_{\ell}(X^{**})$, completely isometrically, via the map $\phi: T\longrightarrow T^{**}$ (see e.g. {\cite[Section 5.3]{BZ2}}). By (ii), $\phi$ is surjective and maps onto the set of complete left $M$-projections. But the left $M$-projections are exactly the contractive projections in $\mathcal{A}_{\ell}(X^{**})$, and since $\mathcal{A}_{\ell}(X^{**})$ is a von Neumann algebra, the span of these projections is dense in it. So $\phi$ maps onto $\mathcal{A}_{\ell}(X^{**})$. 

(v) \ If $X$ is right $M$-embedded then we can show similarly to (iv) that $\mathcal{A}_r(X)\cong \mathcal{A}_r(X^{**})$. By definition, ${Z}(X)=\mathcal{A}_{\ell}(X)\cap \mathcal{A}_r(X)$, hence it follows that ${Z}(X)={Z}(X^{**})$.
\end{proof}

 We will see in {\cite{ABS}} that we can improve (iii) in the above theorem, in certain cases. We will also see that (i) is not true for left $M$-ideals of $X$.

\begin{coro}\label{tensor1}
Let $X$ be a left $M$-embedded operator space. Suppose that $J$ is a complete right $M$-ideal of $X$, and $\otimes_{\beta}$ is any operator space tensor product with the following properties:
\begin{enumerate}
\item[(i)] $-\otimes_{\beta} Id_Z$ is {\em functorial}. That is, if $T:X_1\longrightarrow X_2$ is completely contractive, then $T\otimes_\beta Id_Z: X_1\otimes_\beta Z \longrightarrow X_2\otimes_\beta Z$ is completely contractive,
\item[(ii)] the canonical map $C_2(X)\otimes Z \longrightarrow C_2(X\otimes Z)$extends to a completely isometric isomorphism $C_2(X)\otimes_{\beta} Z \longrightarrow C_2(X\otimes_\beta Z)$,
\item[(iii)] the span of elementary tensors $x\otimes z$ for $x\in X$, $z\in Z$ is dense in $X\otimes_{\beta} Z$. 
\end{enumerate}
Then $J{\otimes_{\beta}}E$ is a complete right $M$-summand of $X{\otimes_{\beta}}E$. In particular, $J\inj{\otimes}E$ is a complete right $M$-summand of $X\inj{\otimes}E$.
\end{coro}
\begin{proof}
Since $J$ is a complete right $M$-ideal of $X$, by Theorem \ref{thm_centralizer}, it is also a right $M$-summand. Hence by the argument in \cite[Section 5.6]{BZ2}, $J{\otimes_{\beta}}E$ is a right $M$-summand of  $X{\otimes_{\beta}}E$. 
\end{proof}

Following is a Banach space result stated for operator spaces. The proof is along similar lines to the Banach space proof. 

\begin{prop}{\label{uep}}
Let $X$ be an operator space and $\pi_{X^*}$ be the canonical projection from $X^{***}$ onto $X^*$. Then the following are equivalent:
\begin{itemize}
\item[(i)] $\pi_{X^*}$ is the only completely contractive projection on $X^{***}$ with kernel $X^{\perp}$. 
\item[(ii)] The only completely contractive operator on $X^{**}$ which restricts to identity on $X$ is $Id_{X^{**}}$.
\item[(iii)] If $U$ is a surjective complete isometry on $X$, then the only completely contractive operator on $X^{**}$ which restricts to $U$ on $X$ is $U^{**}$.
\end{itemize}
\end{prop}

The property in (ii) above is sometimes called the {\em unique extension property}.

\begin{coro}\label{coro_uep}
Every right $M$-embedded space has the unique extension property.
 \end{coro}
\begin{proof}
If $X$ is a right $M$-ideal in $X^{**}$, then by Proposition {\ref{ns}}, $\pi_{X^*}$ is a complete right $L$-projection with kernel $X^{\perp}$. By {\cite[Theorem 3.10(b)] {BEZ}}, it is the only completely contractive projection with kernel $X^{\perp}$. Hence $X$ satisfies all the equivalent conditions in Proposition {\ref{uep}}. In particular, it has the unique extension property. 
\end{proof}

An operator space $X$ has the {\em completely bounded approximation property} (respectively, {\em completely contractive approximation property}) if there exists a net of finite-rank mappings $\phi_{\nu}:X \longrightarrow X$ such that $\norm{\phi_{\nu}}_{cb}\leq K$ for some constant $K$ (respectively, $\norm{\phi_{\nu}}_{cb}\leq 1$) and $\norm{\phi_{\nu}(x)-x}\rightarrow 0$, for every $x\in X$.

\begin{coro}
Let $X$ be a right $M$-embedded operator space. If $X$ has the completely bounded approximation property then $X^*$ has the completely bounded approximation property.
\end{coro}
\begin{proof}
Let $T_{\lambda}$ be a net of finite rank operators in $CB(X)$, such that $\norm{T_{\lambda}}_{cb}\leq K $ for some $K> 0$, and $\norm{T_{\lambda}(x)-x}\longrightarrow 0$. 
We first show that there exists a subnet of $\{T_{\lambda}^*\}$ which converges to $Id_{X^*}$, in the point-weak topology.
We know that $CB(X^{**})$ is a dual operator space with $CB(X^{**})=(X^*\proj{\otimes}X^{**})^*$, so the closed ball of radius $K$ in $CB(X^{**})$, $ K {\rm Ball}(CB(X^{**}))$, is $w^*$-compact. Since $T_{\lambda}^{**}\in K {\rm Ball}(CB(X^{**}))$, there exists a subnet $\{ T_{\lambda_{\nu}}^{**} \}$ and $T$ in  $K {\rm Ball}(CB(X^{**}))$, such that $T_{\lambda_{\nu}}^{**} {\buildrel {w^*} \over \longrightarrow} T $. That is, 
$T_{\lambda_{\nu}}^{**}(\phi)(f) {\longrightarrow} T(\phi)(f)$ for all $\phi\in X^{**}$ and $f\in X^*$.
 In particular for $\hat{x}\in X\subset X^{**}$, the latter convergence implies that $f(T_{\lambda_{\nu}}x)\longrightarrow T(\hat{x})(f)$ for all $f\in X^*$ and $x\in X$. Now since $T_{\lambda} \buildrel \over \longrightarrow Id_X$ in the point-norm topology, it also converges in the point-weak topology. So $f(T_{\lambda_{\nu}}x)  \longrightarrow f(x)$ for all $x\in X$ and $f \in X^*$. Hence $T|_{X}=Id_X$. By Corollary {\ref{coro_uep}}, $X$ has the unique extension property. Hence $T=Id_{X^{**}}$, so $(T_{\lambda_{\nu}}^{**}\phi)(f)\longrightarrow \phi(f)$. Equivalently, $\phi(T_{\lambda_{\nu}}^*f)\longrightarrow \phi(f)$ for all $\phi\in X^{**}$ and $f \in X^*$, which proves the claim. Thus $Id_{X^*}$ is in the point-weak closure of the convex hull of $\{T_{\lambda}^*\}$. But since the norm and the weak topologies coincide on a convex set {\cite[p.477]{DS}}, $Id_{X^*}$ is in the point-norm closure of the convex hull of $\{T_{\lambda}^*\}$.     
\end{proof}

 Along similar lines, we can prove that if a right $M$-embedded space has the completely contractive approximation property then so does its dual. We are grateful to Z. J. Ruan for the following result. Since we could not find this in the literature, we include his proof.

\begin{lemma}{\label{Ruan}}
 Suppose $X^*$ has the completely bounded approximation property and $X$ is a locally reflexive (or $C$-locally reflexive) operator space. Then $X$ has the completely bounded approximation property. 
\end{lemma}
\begin{proof}
We prove the locally reflexive case, the $C$-locally reflexive case is similar.
Suppose that $X$ is locally reflexive. Since $X^*$ has the completely bounded approximation property, there exists a net of finite rank maps $T_{\lambda}: X^*\longrightarrow X^*$ such that $\norm{T_{\lambda}}_{cb}\leq K< \infty$ and $T_{\lambda}\longrightarrow \Id$ in the point-norm topology. Then $\phi_{\lambda}:= (T_{\lambda})^*|_{X}:X\longrightarrow X^{**}$ is a net of finite rank maps such that $\ip{\phi_{\lambda}(x)-x}{f} \rightarrow 0$ for all $x\in X$ and $f\in X^*$. Let $Z_{\lambda}=\phi_{\lambda}(X)$ and $\rho_{\lambda}$ be the inclusion map from $Z_{\lambda}$ to $X^{**}$. Since $\phi_{\lambda}(X)$ is a finite dimensional subspace of $X^{**}$ and $X$ is locally reflexive, for each $\lambda$ we can find a net of completely contractive maps $\rho_{t}^{\lambda}: Z_{\lambda}\longrightarrow X$ such that $\rho_t^{\lambda}$ converges to $\rho_{\lambda}$ in the point-weak$^*$ topology. Then the maps $\psi_{\lambda, t}= \rho_t^{\lambda} \circ \phi_{\lambda}$ are finite rank maps from $X$ to $X$ such that $\norm{\psi_{\lambda, t}}_{cb} \leq K$. Now using a reindexing argument based on {\cite[Lemma 2.1]{B1}}, we show that there exists a net $\gamma$ such that $\lim_{\gamma}\ip{\psi_{\gamma}(x)-x}{f}= 0$ for all $x\in X$ and $f\in X^*$. Define $\Gamma$ to be a set of 4-tuples $(\lambda, t, Y, \epsilon)$, where $Y$ is a finite subset of $X\times X^*$ and where $\epsilon > 0$ is such that $\abs{\psi_{\lambda,t}(x)(f)-\phi_{\lambda}(x)(f)}< \epsilon$ for all $(x,f)\in Y$. Then it is easy to check that $\Gamma$ is a directed set with ordering $(\lambda, t, Y, \epsilon)\leq (\lambda', t', Y', \epsilon')$ iff $\lambda \leq \lambda'$, $Y\subset Y'$ and $\epsilon'\leq \epsilon$. Let $\psi_{\gamma}=\psi_{\lambda,t}$ if $\gamma=(\lambda, t, Y, \epsilon)$. If $\epsilon>0$ choose $\lambda_o$ such that for all $\lambda\geq \lambda_o$ we have $\abs{\phi_{\lambda}(x)(f)-\hat{x}(f)}< \epsilon$. Choose $t_o$ such that $\gamma_o=(\lambda_o, t_o, \{ x,f\}, \epsilon) \in \Gamma$. Now if $\gamma=(\lambda, t, Y, \epsilon')\geq \gamma_o$ then 
\[
\abs{\psi_{\gamma}(x)(f)-\hat{x}(f)} \leq \abs{\psi_{\lambda,t}(x)(f)-\phi_{\lambda}(x)(f)} + \abs{\phi_{\lambda}(x)(f)-\hat{x}(f)} \leq \epsilon' + \epsilon < 2\epsilon. 
\]
Hence $\psi_{\gamma}\rightarrow \Id_X$ in the point-weak topology and thus, $Id_X$ is in the point-weak closure of $K\Ball{CB(X)}$. But the point-weak and the point-norm closures of $K\Ball{CB(X)}$ coincide {\cite[p.477]{DS}}, thus there exist a net $\{\eta_{p}\}\subset K\Ball{CB(X)}$ such that $\eta_{p}\rightarrow \Id _{X}$ in the point-norm topology. 
\end{proof}

{\bf Remark.} A natural question is whether right $M$-embedded or completely $M$-embedded spaces are locally reflexive? Also note that if $X$ has the completely bounded approximation property then by {\cite[Theorem 11.3.3]{ER2}}, $X$ has the strong operator space approximation property. Hence by {\cite[Corollary 11.3.2]{ER2}}, $X$ has the slice map property for subspaces of $B(\ell^{2})$. There seems some hope that the argument in {\cite[Theorem 14.6.6]{ER2}} can be made to imply that $X$ is $1$-exact, and hence is locally reflexive.

\medskip

The following lemma is a well known Banach space result (see {\cite[Lemma III.2.14]{HWW}} for proof). 
\begin{lemma}
For a Banach space $X$ and $x^*\in X^*$ with $\norm{x^*}=1$, the following are equivalent:
\begin{enumerate}
\item[(i)] $x^*$ has a unique norm preserving extension to a functional on $X^{**}$.
\item[(ii)] The relative $w$- and $w^*$-topologies on the ball of $X^*$, $B_{X^*}$ agree at $x^*$, i.e., the map $Id_{B_{X^*}}: (B_{X^*}, w^*)\longrightarrow (B_{X^*}, w)$ is continuous at $x^*$. 
\end{enumerate}
\end{lemma}

\begin{coro}\label{coro_w-w*}
If $X$ is a one-sided $M$-ideal in its bidual, then the relative $w$- and $w^*$-topologies on $B_{X*}$ agree on the unit sphere. 
\end{coro}
\begin{proof}
This is an immediate consequence of the fact that the one-sided $M$-ideals are Hahn-Banach smooth (see e.g.\ {\cite [Chapter 2]{BZ2}}) and the above lemma. 
\end{proof}
The following result follows from Corollary {\ref{coro_w-w*}} (see the argument in {\cite[Corollary III.2.16]{HWW}}). For the definition of a norming subspace see the beginning of this section.

\begin{coro}{\label {norming}}
If $X$ is a one-sided $M$-ideal in its bidual, then $X^*$ contains no proper norming subspace. 
\end{coro}

{\bf Remark.} The above corollary combined with Proposition 2.5 in {\cite{GS}}, immediately gives a second proof of the unique extension property for one-sided $M$-embedded operator spaces. We note that the Proposition 2.5 in {\cite{GS}} is proved for a real Banach space since it uses a lemma ({\cite[Lemma 2.4]{GS}}) on real Banach spaces. However, it is easy to see, using the fact that $(E_{\mathbb{R}})^*=(E^*)_{\mathbb{R}}$, isometrically (see {\cite[Proposition 1.1.6]{Li}}), that the lemma is also true for any complex Banach space $E$. Here $E_{\mathbb{R}}$ denotes the underlying real Banach space.

\begin{prop}{\label{ccc_subspace}}
Let $Y$ be a completely contractively complemented operator space in $Y^{**}$, i.e., $Y\oplus Z =Y^{**}$, and $\norm{[y_{ij}]}\leq \norm{[\phi_{ij}]}$ for all $\phi_{ij}=y_{ij}+{z_{ij}}$ where $y_{ij}\in Y$, $z_{ij}\in Z$ and $\phi_{ij}\in Y^{**}$ for all $i$,$j$. Then $Y$ cannot be a proper right M-ideal in any other operator space.
\end{prop}
\begin{proof}
Let $X$ be an operator space with $Y$ a complete right $M$-ideal in $X$. Suppose that $P$ is a complete left $M$-projection from $ X^{**}$ onto $Y^{\perp\perp}$. By the hypothesis, there is a completely contractive projection $Q:Y^{**}\longrightarrow Y^{**}$ mapping onto $Y$. Let $R$ be the restriction of $(Q\circ P)$ to $X$. Then since $Y^{**}\cong Y^{\perp\perp}$ completely isometrically, $R$ is a completely contractive projection onto Y. Hence by the uniqueness of a left $M$-projection (see e.g. {\cite[Theorem 3.10]{BEZ}}), $R$ has to be a complete left $M$-projection, and thus, $Y$ is a right $M$-summand.
\end{proof}

\begin{prop}
Every non-reflexive right $M$-embedded operator space contains a copy of $c_0$. Moreover, every subspace and every quotient of a right $M$-embedded space, which is not reflexive, contains a copy of $c_0$.
\end{prop}
\begin{proof}
Suppose that $X$ is a non-reflexive right $M$-ideal in its bidual, and suppose that $X$ does not contain a copy of $c_0$. Since $X$ is a $u$-ideal, by {\cite[Theorem 3.5]{GKS}} it is a $u$-summand. Since $u$-summands are contractively complemented, $X$ is the range of a contractive projection. But this implies that $X$ is a right $M$-summand (see the discussion at the end in {\cite[Section 2.3] {BZ2}}). Since $X$ is non-reflexive, it cannot be a non-trivial $M$-summand in $X^{**}$. Hence $X$ has to contain a copy of $c_0$. The rest follows from Theorem {\ref{stability}}.
\end{proof}

Let $X$ be an operator space. Then $\pi_{X^{**}}:=i_{X^{**}}\circ {(i_{X^*})}^*$ is a completely contractive  projection onto $X^{**}$ with kernel $(X^{*})^{\perp}$. So $X^{(4)}=X^{**}\oplus (X^{*})^{\perp}$. The following may be used to give an alternative proof of some results above. 

\begin{prop}{\label{piX}}
If $X$ is a right $M$-embedded operator space, then $\pi_{X^{**}}$ is the only contractive projection from $X^{(4)}$ onto $X^{**}$. 
\end{prop}
\begin{proof}
Since $X$ is right $M$-embedded, then by Theorem {\ref{RNP}}, $X^*$ has the RNP, i.e.,  $(X^*)_{\mathbb{R}}$ has the RNP, where $X_{\mathbb{R}}$ denotes the underlying real Banach space.
 Then by {\cite[p.202 Theorem 3]{DU}}, $\Ball {X^*}_{\mathbb{R}}$ is the closure of the convex hull of its strongly exposed points. If $\psi$ is a strongly exposed point in $\Ball {X^*}_{\mathbb{R}}$, then it is a denting point (see e.g. {\cite{JL}}). Hence $\psi$ is a point of continuity of $\Id : ((X^*)_{\mathbb{R}}, w)\longrightarrow ((X^*)_{\mathbb{R}},\norm{.})$. Thus by {\cite[p.144]{God}}, $(X^{*}_{\mathbb{R}})_{\mathbb{R}}^*$ satisfies the assumptions of \cite[Theorem II.1]{God}. Hence there is a unique contractive $\mathbb{R}$-linear projection from $(X^*_{\mathbb{R}})^{***}$ onto $(X^{*}_{\mathbb{R}})^*$. Since $(X^*)_{\mathbb{R}}=(X_{\mathbb{R}})^*$ ({\cite[Proposition 1.1.6]{Li}}), there is a unique $\mathbb{R}$-linear contractive projection from $(X^{(4)})_{\mathbb{R}}$ to $(X^{**})_{\mathbb{R}}$, and hence a unique $\mathbb{C}$-linear contractive projection from $X^{(4)}$ onto $X^{**}$. 
\end{proof}

{\bf Remark.} Note that the above result also holds for Banach spaces $X$ such that $X^*$ has the RNP, and in particular for $h$-ideals which are strict in the sense of {\cite{GKS}}. It also holds for separable strict $u$-ideals by the proof of Theorem 5.5 from {\cite{GKS}}. 

\section{One-sided $L$-embedded spaces}

We now talk about the dual notion of $L$-embedded spaces. 
For the definition of a one-sided $L$-embedded and a completely $L$-embedded operator space, see Section 2.

\medskip

{\bf Examples.} We list a few examples of right $L$-embedded spaces:
\begin{enumerate}
\item[(a)] Duals of left $M$-embedded spaces.
\item[(b)] Preduals of von Neumann algebras.
\item[(c)] Preduals of subdiagonal operator algebras, in the sense of Arveson {\cite{Arv}}. 
\end{enumerate}

We have already noted (a) in Corollary {\ref{coro_duality}}.
 For (b), note that it is well known that $(\M_*)^{\perp}$ is a $w^*$-closed two-sided ideal in $\M^{**}$, for any von Neumann algebra $M$. So $(\M_*)^{\perp}$ is a complete $M$-ideal in $\M^{**}$. Hence by {\cite[p.8]{BZ2}} and {\cite[Proposition 3.9]{BEZ}}, $\M_*$ is a complete $L$-summand in $\M^*$. 
For (c), let $\A=H^{\infty}(\M,\tau)$, where $\M$ is a von Neumann algebra and $\tau$ a faithful normal tracial state. Then by {\cite[Theorem 2]{Ued}}, $\A$ has a unique predual, namely $\A_*=M_*/\A_{\perp}$. Also, each $M_n(\A)$ is a subdiagonal operator algebra, so applying {\cite[Corollary 2]{Ued}} to $M_n(\A)^*$ we have that each $M_n(\A)_*$ is an $L$-summand in $M_n(\A)^*$. Thus by Lemma {\ref{Lproj}}, $\A_*$ is a complete $L$-summand in $\A^*$.

\medskip

It is shown in \cite{ABS} that if X is right but not left $M$-embedded, then $X^*$ is
left but not right $L$-embedded. Thus the examples mentioned in Section 3
have duals which are left but not right $L$-embedded.


Let $X$ be left $L$-embedded. Then we say a closed subspace $Y$ of $X$ is a {\em left $L$-subspace} if $Y$ is left $L$-embedded and for the right $L$-projection $Q$ from $Y^{**}$ onto $Y$, we have that $\Ker Q=\Ker P\cap Y^{\perp\perp}$, where $P$ is a right $L$-projection from $X^{**}$ onto $X$.

\begin{thm}{\label{subspace_Lembedd}}
Let $X$ be a left $L$-summand in $X^{**}$ and let $Y$ be a subspace of $X$. Let $P:X^{**}\longrightarrow X^{**}$ be a complete right $L$-projection onto $X$. Then the following conditions are equivalent:
\begin{enumerate}
\item[(i)] $Y$ is a left $L$-subspace of $X$.
\item[(ii)] $P(\overline{Y}^{w^*})=Y$.
\item[(iii)] $P(\overline{B_Y}^{w^*})=B_Y$.
\end{enumerate} 
\end{thm}
\begin{proof}
 If $Y$ is a left $L$-subspace then since, $\overline{Y}^{w^*}=Y^{\perp\perp}=Y\oplus (Y^{\perp\perp}\cap \Ker P)$, it is clear that $P(\overline{Y}^{w^*})=Y$. Hence (i) implies (ii). Also since, 
\[
B_Y=P(B_Y)\subset P(\overline{B_Y}^{w^*})=P(B_{Y^{\perp\perp}})=P(B_{\overline{Y}^{w^*}})\subset B_{Y},
\]
it is clear that (ii) and (iii) are equivalent. We now show that (ii) implies (i). Since $P(Y^{\perp\perp})=Y\subset Y^{\perp \perp}$, the restriction of $P$ to $Y^{\perp\perp}=Y^{**}$, say $Q$, is a completely contractive projection from $Y^{**}$ onto $Y$. Also we have that $P\in \mathcal{C}_r(X^{**})$ (for notation see {\cite[Chapter 2]{BZ2}}) and $P=P^{\star}$, and $P(Y^{\perp\perp})\subset Y^{\perp\perp}$, so by {\cite[Corollary 5.12]{BZ2}} we have $Q\in \mathcal{C}_r(Y^{**})$. Thus $Q$ is a right $L$-projection and clearly since $Q=P|_{Y^{\perp\perp}}$, $\Ker Q=\Ker P\cap Y^{\perp\perp}$. Hence $Y$ is a left $L$-subspace of $X$.    
\end{proof}

\begin{coro}
Let $X$ be a left $L$-embedded operator space and $Y$ be a left $L$-subspace of $X$, then $X/Y$ is left $L$-embedded.
\end{coro}
\begin{proof}
Let $P:X^{**}\longrightarrow X^{**}$ be a complete right $L$-projection onto $X$, then by Theorem {\ref{subspace_Lembedd}}, $P$ maps $Y^{\perp\perp}$ onto $Y$. Consider the map
$$P/{Y^{\perp\perp}}:X^{**}/Y^{\perp\perp}\longrightarrow X^{**}/Y^{\perp\perp}$$ given by $(P/{Y^{\perp\perp}})(x^{**}+Y^{\perp\perp})=P(x^{**})+Y^{\perp\perp}$. Then since $P\in \mathcal{C}_r(X^{**})$ (see {\cite[Chapter 2]{BZ2}} for the notation) with $P(Y^{\perp\perp})=P^{\star}(Y^{\perp\perp})\subset Y^{\perp\perp}$, by {\cite [Corollary 5.12]{BZ2}} we have that $P/Y^{\perp\perp}\in \mathcal{C}_r(X^{**}/Y^{\perp\perp})$. So $P/Y^{\perp\perp}$ is a complete right $L$-projection onto $(X+Y^{\perp\perp})/Y^{\perp\perp}$. Since $(X/Y)^{**}$ is completely isometrically isomorphic to $X^{**}/Y^{\perp\perp}$ and under this isomorphism $X/Y$ is mapped onto $(X+Y^{\perp\perp})/Y^{\perp\perp}$, it is clear that $X/Y$ is left $L$-embedded.  
\end{proof}

The following corollary can also be proved using Proposition {\ref{subspace_Lembedd}} (see {\cite[Proposition IV.1.6]{HWW}}).
\begin{coro}
Let $X$ be a left $L$-embedded space and let $Y_1$, $Y_2$, $\{ Y_i \}_{i\in I}$ be left $L$-subspaces of $X$. Then
\begin{enumerate}
\item[(i)] $\cap _{i\in I} Y_i$ is a left $L$-subspace.
\item[(ii)] $Y_1+Y_2$ is closed if and only if $Y_1+Y_2$ is a left $L$-subspace of $X$.
\end{enumerate} 
\end{coro}

We omit the proofs of the proposition below, because it is identical to the classical version (see {\cite[Proposition IV.1.12]{HWW}}).

\begin{prop}
Let $X$ be a left $L$-embedded space and let $Y$ be a left $L$-subspace of $X$. Then $Y$ is proximinal in $X$ and the set of best approximations to $x$ from $Y$ is weakly compact for all $x$ in $X$.
\end{prop}

\medskip
The following two results are non-commutative versions of some of Godefroy's results.
\begin{prop}
Let $X$ be left $L$-embedded and let $P$ be a complete right $L$-projection from $X^{**}$ onto $X$. Then
\begin{enumerate}
\item[(i)] there is at most one predual of $X$, up to complete isometric isomorphism, which is right $M$-embedded,
\item[(ii)] there is a predual of $X$ which is a right $M$-ideal in its bidual if and only if $\Ker{P}$ is $w^*$-closed in $X^{**}$.
\end{enumerate}
\end{prop}
\begin{proof}
(i) \ Let $Y_1$ and $Y_2$ be two preduals of $X$, that is $Y_1^*\cong X\cong Y_2^*$ completely isometrically via a map $I:Y_1^*\longrightarrow Y_2^*$. Let $P={I^{**}}^{-1} \pi_{{Y_2}^*} I^{**}$, then $P:Y_1^{***}\longrightarrow Y_1^{***}$ is a completely contractive projection onto $Y_1^{*}$. Thus by {\cite[Theorem 3.10(a)]{BEZ}}  $\pi_{{Y_1}^*}=P={I^{**}}^{-1} \pi_{{Y_2}^*} I^{**}$. By basic functional analysis, this is equivalent to the $w^*$-continuity of $I$, which implies that $I=J^*$ for some complete isometric isomorphism $J:Y_2\longrightarrow Y_1$. Thus the predual is unique up to complete isometry.

(ii) \ Suppose that $Y$ is a right $M$-embedded operator space such $Y^*=X$. Then $\pi_{Y^*}$ is a complete right $L$-projection from $Y^{***}$ onto $Y^*=X$, so $\pi_{Y^*}=P$. The kernel of $\pi_{Y^*}$ is $i_Y(Y)^{\perp}\subset Y^{***}$, which is clearly $w^*$-closed in $X^{**}$. Conversely, suppose that $\Ker P$ is $w^*$-closed in $X^{**}$. Let $Y=(\Ker P)_{\perp}\subset X^*$, then $Y^{\perp}=\overline{\Ker P}^{w^*}=\Ker P=$ Ran$(I-P)$, which means that $Y^{\perp}$ is a left $L$-summand. Hence $Y$ is a right $M$-ideal in $X^*$. Since $P:X^{**}\longrightarrow X^{**}$ is a complete quotient map onto $X$, $X^{**}/{\Ker P}\cong X$. But $X^{**}/{\Ker P}\cong ((\Ker P)_{\perp})^*$, so $X\cong Y^*$.  
\end{proof}

\begin{coro} Let $X$ be a right $M$-ideal in its bidual and $Y$ be a $w^*$-closed subspace of $X^*$. Then
\begin{enumerate}
\item[(i)] $Y$ is the dual of a space which is a right $M$-ideal in its bidual.
\item[(ii)] $Y$ is a left $L$-summand in its bidual.
\end{enumerate}
\end{coro}

\begin{proof}
  If $Y$ is $w^*$-closed, then $Y=(X/Y_{\perp})^*$. Now since $X$ is right $M$-embedded, by Theorem {\ref{stability}}, $X/Y_{\perp}$ is right $M$-embedded. This proves (i).    
It is easy to see that (ii) follows by Corollary {\ref {coro_duality}}.
\end{proof}

\medskip

Analogues of many classical results about the RNP are also true for the right $L$-embedded spaces. For instance, suppose that $X$ is right $L$-embedded and $P$ is a left $L$-projection from $X^{**}$ onto $X$. Then, if the ball of $\Ker P$ is $w^*$-dense in the ball of $X^{**}$, $X$ fails to have the RNP. This is because the unit ball of $X$ does not have any strongly exposed points (see e.g.\ {\cite[Remark IV.2.10 (a)]{HWW}}).

\begin{prop}
Let $X$ be a left $L$-embedded operator space and $Y\subset X$ be a left $L$-subspace. Let $Z$ be an operator space such that $Z^*$ is an injective Banach space $($resp.\ injective operator space$)$. Then for every contractive $($resp.\ completely contractive$)$ operator $T:Z\longrightarrow X/Y$ there exists a contractive $($resp.\ completely contractive$)$ map $S: Z\longrightarrow X$ such that $qS=T$, where $q :X\longrightarrow X/Y$ is the canonical quotient map.
\end{prop}
The above proposition can be proved by routine modifications to the argument in {\cite[Proposition IV.2.12]{HWW}}. For the following corollaries see arguments in {\cite[Corollary IV.2.13]{HWW}} and {\cite[Corollary IV.2.14]{HWW}}, respectively.

\begin{coro}
If $X$ is a right $L$-embedded space with $Y$ a left $L$-subspace of $X$, and if $X/Y$ contains a subspace $W$ isometric $($resp.\ completely isometric$)$ to $L^1(\mu)$, then there is a subspace $Z$ of $X$ such that $q(Z)=W$ and $q|_{Z}$ is an isometric $($resp.\ completely isometric$)$ isomorphism. If also $X/Y\cong L^1(\mu)$, then there is a contractive $($resp.\ completely contractive$)$ projection $P$ on $X$ with $Y=\Ker P$.
\end{coro}

\begin{coro}
Let $X$ be a right $L$-embedded space with $Y$ a left $L$-subspace of $X$. Then if $X$ has the RNP then $X/Y$ has the RNP.
\end{coro}

{\bf Acknowledgments.} This work was done as a part of the Ph.D. thesis of the author at the University of Houston. We thank our Ph.D. adviser, Dr.\ David Blecher, for proposing this project and continually supporting the work. We are grateful for his insightful comments and very many suggestions and corrections.

\end{document}